How to Construct a Large Table of Reciprocals of Babylonian Mathematics

Kazuo MUROI

§1. Introduction

A large table of reciprocals of the Seleucid period (358 B.C.–63 B.C.), AO 6456, was published in 1922 by F. Thureau-Dangin through his beautiful hand-copy of the clay tablet.[1] Ten years after, the first fairly reliable study of the table was presented by O. Neugebauer, and almost all the mistakes of the numbers of the tablet, which amount to 56 in number, were corrected by him.[2] Thanks to Neugebauer's study, we can mathematically understand the large table of reciprocals of AO 6456 except for one point: how did the Babylonian scribe of our tablet construct such a large table?

In the present paper I shall clarify the Babylonian method by which many large tables of reciprocals could be constructed.

§2. The outline of AO 6456

The tablet AO 6456 begins with a prayer:

*ina a-mat $^{d}$Anu u An-tum mim-ma ma-la dù-uš ina* šu$^{II}$*-ia liš-lim*

"By means of a command of (the gods) Anu and Antum, everything that I do should be perfect in my hands."

After this set phrase, which is characteristic of Uruk texts, 157 pairs of a number and its reciprocal are listed, the numbers running from 1 to 3:

igi-1-gál-bi 1-àm "The reciprocal of 1 is 1."

igi 1;0,16,53,53,20     0;59,43,10,50,52,48

igi 1;0,40,53,20        0;59,19,34,13,7,30

igi 1;0,45              0;59,15,33,20

| | |
|---|---|
| igi 1;1,2,6,33,45 | 0;58,58,56,38,24 |
| igi 1;1,26,24 | 0;58,35,37,30 |
| igi 1;1,30,33,45 | 0;58,31,39,35,18,31,6,40 |
| … … | … … |
| igi 1;29,12,19,26,34,23,19,49,38,8,36,52,20,44,26,40 | 0;40,21,22,41,0,9 |
| … … | … … |
| igi 2;15 | 0;26,40 |
| … … | … … |
| igi 2;55,46,52,30 | 0;20,28,48 |
| igi 2;57,46,40 | 0;20,15 |

Under the catchline of igi-3-gál-bi 0;20, the text ends with the following colophon:

> *per-su reš-tu-u*∶1∶*a-mu-ú*∶2∶*a-mu-ú* nu al-til / im ᵐ*Nidintu-Anu* a *šá* ᵐ*Ina-qí-bit-Anu* a ᵐ*Hun-zu-u* ˡᵘ*maš-maš Anu u An-tum* Uruk^(ki)-*u* / qàt ᵐ*Ina-qí-bit-Anu* dumu-a-ni-*šú*

> "The first section. '1' is a head number. '2' is a head number. It is not completed. / The tablet of Nidintu-Anu who is a son of Inaqibit-Anu (who is) a son of Hunzû, the incantation priest of (the gods) Anu and Antum in Uruk. / (By) the hand of Inaqibit-Anu, his son, (it was written)."[3]

Thus the scribe of this tablet, Inaqibit-Anu, calculated even the reciprocal of a sexagesimal seventeen-place number. It is highly improbable that he calculated these reciprocal pairs one by one, and there must have been a systematic method from which 157 pairs of numbers were derived.

§3. Successive doubling or successive tripling

Recently I have investigated three mathematical technical terms, a-rá-kár, a-rá-hi, and a-rá-gub-ba, and I have clarified that each of these terms played an important role in Babylonian number theory.[4] Through the terms, especially a-rá-kár and a-rá-gub-ba, we can also find an answer to our question; how did the Babylonians construct such a large table of reciprocals as AO 6456?

The Babylonians occasionally called the positive integers 2, 3, and 5 a-rá-gub-ba "normal factor" because each reciprocal of the numbers can be expressed by a finite sexagesimal fraction. In other words, they knew the fact that the reciprocal of a number $n$ is a finite sexagesimal number if and only if $n = 2^{\alpha} \cdot 3^{\beta} \cdot 5^{\gamma}$ ($\alpha, \beta, \gamma$ : integers). Moreover, they deliberately used the numbers 2 and 3 in order to construct large numbers, that is to say, the Babylonians obtained large numbers by successive doubling or successive tripling which was called a-rá-kár "blowing up multiplication". The opposite of this a-rá-kár operation is probably igi-te-en "deflation of foam"[5] which usually means "fraction, proportion" both in mathematical texts and non-mathematical texts. By the successive doubling (a-rá-kár) and the successive halving (igi-te-en) the Babylonians were able to make many tables of reciprocals, for example:

$m$  igi-bi $\overline{m}$ , "$m$ , its reciprocal is $\overline{m}$ ."

$m = 2,5 \cdot 2^n$  (n=1  2  3  ...  ...  30).

Out of thirty pairs of these numbers we can select ten pairs that also occur in AO 6456 if we adjust the sexagesimal point properly:

$\overline{1;48} = 0;33,20$  from  $\overline{33,20} = 0;0,1,48$  (n=4)

$\overline{1;6,40} = 0;54$ from $\overline{1,6,40} = 0;0,0,54$ (n=5)

$\overline{2;13,20} = 0;27$ from $\overline{2,13,20} = 0;0,0,27$ (n=6)

$\overline{1;41,15} = 0;35,33,20$ from $\overline{35,33,20} = 0;0,0,1,41,15$ (n=10)

$\overline{1;11,6,40} = 0;50,37,30$ from $\overline{1,11,6,40} = 0;0,0,0,50,37,30$ (n=11)

$\overline{1;34,55,18,45} = 0;37,55,33,20$ from $\overline{37,55,33,20} = 0;0,0,0,1,34,55,18,45$ (n=16)

$\overline{1;15,51,6,40} = 0;47,27,39,22,30$ from $\overline{1;15,51,6,40}$

$\qquad = 0;0,0,0,0,47,27,39,22,30$ (n=17)

$\overline{2;31,42,13,20} = 0;23,43,49,41,15$ from $\overline{2,31,42,13,20}$

$\qquad = 0;0,0,0,0,23,43,49,41,15$ (n=18)

$\overline{1;28,59,21,19,41,15} = 0;40,27,15,33,20$ from $\overline{40,27,15,33,20}$

$\qquad = 0;0,0,0,0,1,28,59,21,19,41,15$ (n=22)

$\overline{1;20,54,31,6,40} = 0;44,29,40,39,50,37,30$ from $\overline{1,20,54,31,6,40}$

$\qquad = 0;0,0,0,0,0,44,29,40,39,50,37,30$ (n=23). See Table 1.

Table 1

| n | $2^n \cdot 5^3$ | $\overline{2^n \cdot 5^3}$ |
|---|---|---|
| 1 | 4,10 | 0;0,14,24 |
| 2 | 8,20 | 0;0,7,12 |
| 3 | 16,40 | 0;0,3,36 |
| 4 | 33,20 | 0;0,1,48 |
| 5 | 1,6,40 | 0;0,0,54 |
| 6 | 2,13,20 | 0;0,0,27 |
| 7 | 4,26,40 | 0;0,0,13,30 |

| 8  | 8,53,20              | 0;0,0,6,45                          |
|----|----------------------|-------------------------------------|
| 9  | 17,46,40             | 0;0,0,3,22,30                       |
| 10 | 35,33,20             | 0;0,0,1,41,15                       |
| 11 | 1,11,6,40            | 0;0,0,0,50,37,30                    |
| 12 | 2,22,13,20           | 0;0,0,0,25,18,45                    |
| 13 | 4,44,26,40           | 0;0,0,0,12,39,22,30                 |
| 14 | 9,28,53,20           | 0;0,0,0,6,19,41,15                  |
| 15 | 18,57,46,40          | 0;0,0,0,3,9,50,37,30                |
| 16 | 37,55,33,20          | 0;0,0,0,1,34,55,18,45               |
| 17 | 1,15,51,6,40         | 0;0,0,0,0,47,27,39,22,30            |
| 18 | 2,31,42,13,20        | 0;0,0,0,0,23,43,49,41,15            |
| 19 | 5,3,24,26,40         | 0;0,0,0,0,11,51,54,50,37,30         |
| 20 | 10,6,48,53,20        | 0;0,0,0,0,5,55,57,25,18,45          |
| 21 | 20,13,37,46,40       | 0;0,0,0,0,2,57,58,42,39,22,30       |
| 22 | 40,27,15,33,20       | 0;0,0,0,0,1,28,59,21,19,41,15       |
| 23 | 1,20,54,31,6,40      | 0;0,0,0,0,0,44,29,40,39,50,37,30    |
| 24 | 2,41,49,2,13,20      | 0;0,0,0,0,0,22,14,50,19,55,18,45    |
| 25 | 5,23,38,4,26,40      | 0;0,0,0,0,0,11,7,25,9,57,39,22,30   |
| 26 | 10,47,16,8,53,20     | 0;0,0,0,0,0,5,33,42,34,58,49,41,15  |
| 27 | 21,34,32,17,46,40    | 0;0,0,0,0,0,2,46,51,17,29,24,50,37,30 |
| 28 | 43,9,4,35,33,20      | 0;0,0,0,0,0,1,23,25,38,44,42,25,18,45 |
| 29 | 1,26,18,9,11,6,40    | 0;0,0,0,0,0,0,41,42,49,22,21,12,39,22,30 |
| 30 | 2,52,36,18,22,13,20 | 0;0,0,0,0,0,0,20,51,24,41,10,36,19,41,15 |

In this way all the reciprocal pairs of AO 6456 can be obtained from a set of basic tables of reciprocals like the above one.  By my count, the following twenty-three tables will suffice to reproduce the table of AO 6456:

$(2^n, 2^{-n})$   $(2^n \cdot 3, 2^{-n} \cdot 3^{-1})$   $(2^n \cdot 3^2, 2^{-n} \cdot 3^{-2})$ … …  $(2^n \cdot 3^8, 2^{-n} \cdot 3^{-8})$   $(2^n \cdot 5, 2^{-n} \cdot 5^{-1})$  …

$(2^n \cdot 5^{12}, 2^{-n} \cdot 5^{-12})$ $(3^n, 3^{-n})$ $(3^n \cdot 5, 3^{-n} \cdot 5^{-1})$ for n=1   2   3   …   …   30.

If we pick up the 157 reciprocal pairs in ascending order from the above 23 basic reciprocal tables, we can obtain the same table as that of AO 6456.  Moreover, if we calculate these numbers by hand through the successive doubling (tripling) and the successive halving (division by 3), we will realize that some seemingly tremendous reciprocal pairs of AO 6456 are easily obtained, and confirm that the *blowing up multiplication* (a-rá-kár) with the *deflation of foam* (igi-te-en) also played an important role in constructing large tables of reciprocals.

Notes

(1) F. Thureau-Dangin, *Tablettes d'Uruk,* 1922, Plates 55-58.

(2) O. Neugebauer, Sexagesimalsystem und babylonische Bruchrechnung IV, *Quellen und Studien zur Geschichte der Mathematik Astronomie und Physik*, Abteilung B; Studien, Band 2, 1932-1933, pp. 199-210.

　　O. Neugebauer, *Vorgriechische Mathematik*, 1934, pp. 4-16.

　　O. Neugebauer, *Mathematische Keilschrift-Texte*, Erster Teil, 1935, pp. 14-22.

(3) K. Muroi, The Colophon in the Table of Reciprocals of the Seleucid Period (written in Japanese), *Kagakusi-Kenkyu* vol. 32, 1993, pp. 29-31.

(4) K. Muroi, *Studies in Babylonian Mathematics*, No. 3, 2007, pp. 2-8.

(5) For the meaning of "bubble, foam" of igi, see; A. L. Oppenheim and others, *The Assyrian Dictionary* (=CAD), vol. 7　I /J, 1960, pp. 153-158 (*īnu*).　Also, CAD, vol. 6, H, 1956, p. 250 (*hurhummatu*).　For the meaning of "to burst (said of bubbles)" of te-en, see CAD, vol. 2, B, 1965, pp. 72-74 (*balû*).